# Optimal Design and Control of 4-IWD Electric Vehicles Based on a 14-DOF Vehicle Model

Huilong Yu , *Member, IEEE*, Federico Cheli, and Francesco Castelli-Dezza , *Member, IEEE*

*Abstract*—A 4-independent wheel driving (4-IWD) electric vehicle has distinctive advantages with both enhanced dynamic and energy efficiency performances since this configuration provides more flexibilities from both the design and control aspects. However, it is difficult to achieve the optimal performances of a 4-IWD electric vehicle with conventional design and control approaches. This paper is dedicated to investigating the vehicular optimal design and control approaches, with a 4-IWD electric race car aiming at minimizing the lap time on a given circuit as a case study. A 14-DOF vehicle model that can fully evaluate the influences of the unsprung mass is developed based on Lagrangian dynamics. The 14-DOF vehicle model implemented with the reprogrammed Magic Formula tire model and a time-efficient suspension model supports metric operations and parallel computing, which can dramatically improve the computational efficiency. The optimal design and control problems with design parameters of the motor, transmission, mass center, anti-roll bar and the suspension of the race car are successively formulated. The formulated problems are subsequently solved by directly transcribing the original problems into large-scale nonlinear optimization problems based on trapezoidal approach. The influences of the mounting positions of the propulsion system, the mass and inertia of the unsprung masses, the anti-roll bars, and suspensions on the lap time are analyzed and compared quantitatively for the first time. Some interesting findings that are different from the 'already known facts' are presented.

## I. INTRODUCTION

DEVELOPING electric vehicles (EVs) has been globally recognized as a promising solution to face the challenges of air pollution, fossil oil crisis, and greenhouse gas emissions, which leads to mushroomed penetration of EVs in the last decade [1]. Most of the research interests on EVs are focused on energy management [2], [3], electric motor control [4]–[6] and dynamic control [7]–[9]. Continuous research on these topics have improved the energy efficiency and dynamic performance of the EVs a lot. However, optimal design and control of EVs that can help to improve the performances further are seldom discussed in existing literature. With the fast advances in electric vehicles and the emerging autonomous electric vehicles, optimal control theory will undoubtedly play more important role in realizing a 'one line one design, one line one control' concept with assist of the continuously reduced manufacturing cost, in which condition the optimality of the concerned performance is more meaningful.

Normally, the proper parameters of the motor and transmission are the primary considerations to meet the performance requirements in electric vehicle design [10]. This is due to the fact that the power density, dynamic performance, energy efficiency and cost of an electric powertrain rely heavily on the matching of the motor and transmission. However, most of the existing work designed the electric powertrain following a conventional method which is unlikely to obtain the optimal performances [11]. This common approach can be summarized into three steps. The first step is to define the motor power according to the requirements on dynamic performance with a simplified point-mass vehicle model. Then the motor will be chosen from the available products of the manufacturers considering the power density, cost, etc. The last step is to select the gearbox according to the torque-speed characteristics of the motor, the required maximum speed and the maximum torque on the wheels. The limitations of the conventional approach are: (a) the most frequently employed simplified point-mass model can not predict the vehicle dynamic performances more precisely, e.g., influences of the mounting positions of the electric motors can not be evaluated, and the number of parameters can be optimized is limited; (b) the optimal design solutions can not be obtained with this manually design method. For a 4-IWD electric vehicle, the motors and transmissions can be put onboard or in-wheel. However, the influences of their mounting positions on the lap time have not been evaluated quantitatively in the existing work. In addition, the chassis design which is known to play a significant role in the vehicle dynamic performances [12]–[14], nonetheless, is mostly designed following the engineering experience based on manual calculation [15]. There is a certain amount of existing work concerned the optimal control problems of vehicles [16]–[18], however, vehicular optimal design is seldom discussed. In model based design and control methodology, the modelling work serving as the base is particularly of great significance, however, the mostly

Manuscript received March 25, 2018; revised July 12, 2018; accepted September 2, 2018. Date of publication September 17, 2018; date of current version November 12, 2018. The work of H. Yu is supported by the China Scholarship Council. The review of this paper was coordinated by Dr. D. Cao. *(Corresponding author: Federico Cheli.)*

H. Yu was with the Department of Mechanical Engineering, Politecnico di Milano, Milano 20156, Italy. He is now with Vehicle Intelligence Pioneers Inc., and also with Qingdao Academy of Intelligent Industries, Qingdao 266109, P.R. China (e-mail: huilong.yu@ieee.org).

F. Cheli and F. Castelli Dezza are with the Department of Mechanical Engineering, Politecnico di Milano, 20156, Milano, Italy (e-mail: federico.cheli@polimi.it; francesco.castellidezza@polimi.it).

Color versions of one or more of the figures in this paper are available online at http://ieeexplore.ieee.org.

Digital Object Identifier 10.1109/TVT.2018.2870673

implemented 2-DOF, 3-DOF and 7-DOF vehicle models can not represent the practical vehicle behavior precisely and it is only possible to optimize limited number of parameters due to their simplifications though they can save a lot of computing efforts [19].

This work aims to propose an optimal design approach to overcome the aforementioned drawbacks based on a developed 14-DOF vehicle model with improved computing efficiency. In particular, the design of a 4-IWD electric race car aiming at minimizing the lap time on a given circuit is investigated as a case study. In order to test the performance of a design solution, a corresponding control strategy should be developed. However, there are various kinds of control approaches for electric vehicles [20]–[22] and accordingly, different control strategies may result in different results even with the same designed race car. Thus, the optimal control of the electric race car is coupled into the optimal design problem in this work which is reasonable in practice. The final results will include both the optimal design and control solutions.

There novelty and original contributions of this work with respect to the existing literature are presented as followings: First, a vectorized 14-DOF vehicle model that can fully evaluate the influences of the unsprung mass is developed in MATLAB based on Lagrangian dynamics. In particular, the 14-DOF vehicle model is implemented with the reprogrammed full set Magic Formula tire model [23] that supports '.tir' tire data file as input and metric operations to improve the computation efficiency. A time-efficient suspension model is also developed to describe the relationships between the wheel jounce and spring force, damping force, toe angle, steering angle, camber angle, etc. Second, the optimal design and control problems with parameters of the propulsion system, the mass center, the anti-roll bar and the suspension of the electric race car as design parameters are successively formulated in standard formats based on the developed 14-DOF vehicle model and a path following model in curvilinear coordinate system for the first time. Third, the complicated large-scale optimal design and control problems based on the 14-DOF vehicle model are solved based on direct transcription methods for the first time with respect to the existing efforts. Fourth, results of different optimization cases and the influences of the mounting positions of the propulsion system, the mass and inertia variation of unsprung masses on the lap time are analyzed and compared quantitatively, which is rarely found in the search-able literature. Some new findings that are different from the facts that addressed by most engineers and researchers are presented.

The remainder of this work is organized as follows. Section II details the formulation of the optimal design and control problem, with the objective, variables and constraints are presented. Section III elaborates the derivation and validation of the 14-DOF vehicle model in different maneuvers. Section IV depicts briefly the employed numerical optimal control approach, the simulation parameters and the optimization settings. Section V and Section VI demonstrated the obtained optimal parameters, coupled with trajectory, control and state variables. Results of different optimization cases are compared and analyzed. Section VII concludes this work.

## II. PROBLEM FORMULATION

This section gives an overall description of the optimal design and control problem. In this work, the objective is to minimize the lap time $t_f$:

$$J = \min t_f \tag{1}$$

subject to:
- the first order dynamic constraints

$$\dot{\boldsymbol{x}}(t) = \boldsymbol{f}[\boldsymbol{x}(t), \boldsymbol{u}(t), t, \boldsymbol{p}] \tag{2}$$

- the boundaries of the state, control and design variables

$$\boldsymbol{x}_{\min} \leqslant \boldsymbol{x}(t) \leqslant \boldsymbol{x}_{\max}$$
$$\boldsymbol{u}_{\min} \leqslant \boldsymbol{u}(t) \leqslant \boldsymbol{u}_{\max}$$
$$\boldsymbol{p}_{\min} \leqslant \boldsymbol{p} \leqslant \boldsymbol{p}_{\max} \tag{3}$$

- the algebraic path constraints

$$\boldsymbol{g}_{\min} \leqslant \boldsymbol{g}[\boldsymbol{x}(t), \boldsymbol{u}(t), t, \boldsymbol{p}] \leqslant \boldsymbol{g}_{\max} \tag{4}$$

- and the boundary conditions:

$$\boldsymbol{b}_{\min} \leqslant \boldsymbol{b}[\boldsymbol{x}(t_0), t_0, \boldsymbol{x}(t_f), t_f, \boldsymbol{p}] \leqslant \boldsymbol{b}_{\max} \tag{5}$$

where $\dot{\boldsymbol{x}}$ is the first order derivative of the state variables, $\boldsymbol{f}$ is the dynamic model, $\boldsymbol{x}, \boldsymbol{u}, \boldsymbol{p}$ are respectively the state, control and design vector with their lower and upper bounds: $\boldsymbol{x}_{\min}, \boldsymbol{u}_{\min}, \boldsymbol{p}_{\min}$ and $\boldsymbol{x}_{\max}, \boldsymbol{u}_{\max}, \boldsymbol{p}_{\max}$. While $\boldsymbol{g}$ and $\boldsymbol{b}$ are the path and boundary equations respectively with their lower and upper bounds $\boldsymbol{g}_{\min}, \boldsymbol{b}_{\min}$ and $\boldsymbol{g}_{\max}, \boldsymbol{b}_{\max}$. The dimensions of the input and output variables in Equations (2), (4) and (5) are separately given as:

$$\boldsymbol{f}: \mathbb{R}^{n_x} \times \mathbb{R}^{n_u} \times \mathbb{R} \times \mathbb{R}^{n_p} \to \mathbb{R}^{n_x}$$
$$\boldsymbol{g}: \mathbb{R}^{n_x} \times \mathbb{R}^{n_u} \times \mathbb{R} \times \mathbb{R}^{n_p} \to \mathbb{R}^{n_g}$$
$$\boldsymbol{b}: \mathbb{R}^{n_x} \times \mathbb{R} \times \mathbb{R}^{n_x} \times \mathbb{R} \times \mathbb{R}^{n_p} \to \mathbb{R}^{n_b} \tag{6}$$

The state variables $\boldsymbol{x}$, control variables $\boldsymbol{u}$, design parameters $\boldsymbol{p}$ are described in the following paragraphs.

### A. Variables

The state vector $\boldsymbol{x}$ includes the 14 DOF and their derivatives of the vehicle model, and 3 additional variables to describe the vehicle position in curvilinear coordinate system, so the number of the state variables $n_x = 31$ and $\boldsymbol{x}$ is denoted as:

$$\begin{aligned}\boldsymbol{x} = \{&\dot{X}_{A,b}, \dot{Y}_{A,b}, \dot{Z}_{A,b}, \dot{\varphi}, \dot{\phi}, \dot{\psi}, \dot{z}_{fr}, \dot{z}_{fl}, \dot{z}_{rr}, \dot{z}_{rl},\\ &\dot{\theta}_{fr}, \dot{\theta}_{fl}, \dot{\theta}_{rr}, \dot{\theta}_{rl}, X_{A,b}, Y_{A,b}, Z_{A,b}, \varphi, \phi, \psi,\\ &z_{fr}, z_{fl}, z_{rr}, z_{rl}, \theta_{fr}, \theta_{fl}, \theta_{rr}, \theta_{rl}, s, n, \chi\}\end{aligned} \tag{7}$$

where $X_{A,b}$, $Y_{A,b}$ and $Z_{A,b}$ are the displacements of the mass center in longitudinal, lateral and vertical directions of the global reference system, $\varphi$, $\phi$, and $\psi$ are the roll angle, pitch angle and yaw angle of the vehicle body, $z_i$ and $\theta_i$ are the vertical and rotational displacement of each wheel, $s$, $n$, $\chi$ are the traveled distance, normal distance to the reference trajectory and orientation angle of the vehicle in the curvilinear coordinate system, respectively. In this work, $i = \{fr, fl, rr, rl\}$ means front right, front left, rear right and rear left.

The race car is assumed to be controlled with front wheel steering and four independent wheel driving. The dynamic responses of the steering system and the electric motors are not taken into account in this study, the corresponding control vector is:

$$\boldsymbol{u} = [\delta, T_{fr}, T_{fl}, T_{rr}, T_{rl}], \quad n_u = 5 \quad (8)$$

where $\delta$ is the steering angle, $T_i$ is the driving/braking torque acted on each wheel.

The design variables of the 4-IWD electric race car are presented respectively in Section V and VI according to different optimization cases.

### B. Algebraic Path Constraints

The algebraic path constraints are a set of constraints that can be denoted as functions of the state, control, final time and design parameters. For the motor design, the maximum rotational speed $N_{\max,i}$ should be constrained to a user set range $[\boldsymbol{cl}_{N_{\max}}, \boldsymbol{cu}_{N_{\max}}]$,

$$\boldsymbol{cl}_{N_{\max}} \leq \boldsymbol{N}_{\max} = \boldsymbol{N}_b \boldsymbol{\beta} \leq \boldsymbol{cu}_{N_{\max}} \quad (9)$$

where $N_b$ is the base speed and $\beta$ is the constant power speed ratio (CPSR) of each motor.

In order to let the motors work within their available operation zone, the motor speed $\boldsymbol{N}_m$ and output torque $\boldsymbol{T}_m$ should be respectively constrained within the user set lower and upper bounds:

$$\boldsymbol{cl}_{N_{cm}} \leq \boldsymbol{N}_{cm} = \boldsymbol{N}_m - \boldsymbol{N}_{\max} \leq \boldsymbol{cu}_{N_{cm}}$$
$$\boldsymbol{cl}_{T_{cm}} \leq \boldsymbol{T}_{cm} = \boldsymbol{T}_m - \boldsymbol{T}_{\max} \leq \boldsymbol{cu}_{T_{cm}} \quad (10)$$

where $N_{cm}$ and $T_{cm}$ are respectively the constraint function related with the motor speed and torque, the units of torque, rotation speed, power are respectively $Nm$, $rpm$ and $kW$, the available maximum torque of the motor is given as a function of the gear ratio $i_g$, maximum power $P_{\max}$ and the motor speed:

$$\boldsymbol{T}_{\max} = \begin{cases} \dfrac{9550 \boldsymbol{i}_g \boldsymbol{P}_{\max}}{\boldsymbol{N}_b}, & \boldsymbol{N}_m \leq \boldsymbol{N}_b \\ \dfrac{9550 \boldsymbol{i}_g \boldsymbol{P}_{\max}}{\boldsymbol{N}_m}, & \boldsymbol{N}_m > \boldsymbol{N}_b \end{cases} \quad (11)$$

The normal load $\boldsymbol{F}_z$, the tire slip $\boldsymbol{\kappa}$, $\boldsymbol{\alpha}$ and also should be constrained within their allowable range:

$$\boldsymbol{ccl}_{F_z} \leq \boldsymbol{F}_z \leq \boldsymbol{cu}_{F_z}$$
$$\boldsymbol{cl}_{\kappa} \leq \boldsymbol{\kappa} \leq \boldsymbol{cu}_{\kappa}$$
$$\boldsymbol{cl}_{\alpha} \leq \boldsymbol{\alpha} \leq \boldsymbol{cu}_{\alpha} \quad (12)$$

where $\boldsymbol{cl}_i$ and $\boldsymbol{cu}_i$ means the minimum and maximum value of the mentioned variables, respectively.

### III. VEHICLE MODELLING

The configuration of the entire vehicle model is presented in Fig. 1, where the interactions between the vehicle body,

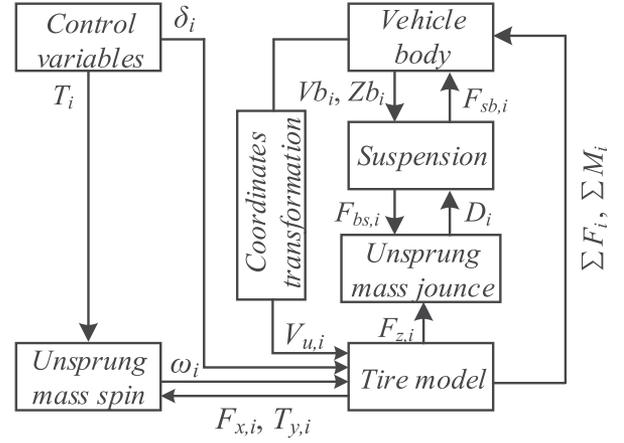

Fig. 1. Vehicle model configuration.

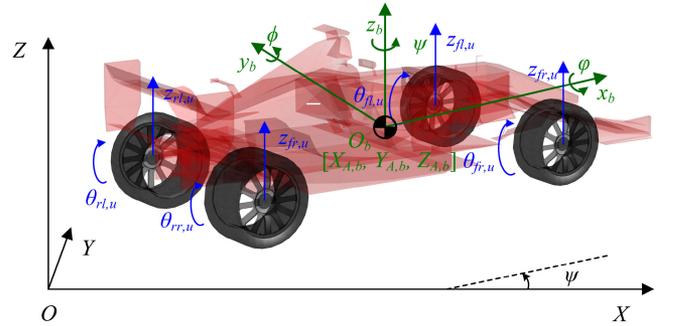

Fig. 2. Degree of freedoms of the 14-DOF vehicle model.

suspension, unsprung mass and tire model are demonstrated. The inputs of the developed 14-DOF vehicle model are the steer wheel angle $\delta_i$ and wheel torque $T_i$. The spin motion of each wheel is driven by the torque and longitudinal force acted on each. The inputs of each tire model are the angular velocity of the wheel $\omega_i$, wheel center velocity $V_{u,i}$, camber angle and normal load $F_{z,i}$, while the outputs are the tire forces and moments. The vertical motion of the each unsprung mass is driven by the normal force of the tire road interaction and the vertical suspension force. While the 6 DOF of vehicle body are driven by the longitudinal, lateral tire forces, aerodynamic forces and the suspension forces acted on it. The vehicle model are derived according to the following assumptions: 1) The inertial matrix of the body is assumed to be diagonal in the reference system fixed to the body and its origin is the center of the gravity; 2) The roll and pitch angles of the chassis are supposed to be small enough to be considered independent from the other; 3) The kinetic energy of both camber rotation and steering angle is neglected; 4) The gyroscopic effect of each wheel is neglected.

### A. 14-DOF Vehicle Model

The 14-DOF vehicle model in this work represents the dynamic behavior of a simplified vehicle consists of five rigid parts, one of which is the vehicle body (sprung mass), and the left four are the connected four wheel parts (unsprung mass). As it is shown in Fig. 2, there are 6 DOF of the vehicle body allows it to displace in the longitudinal, lateral and vertical direction as weel as to roll, pitch and yaw. In this work, the four wheels are supposed to be fixed with the chassis to move in the longitudinal and lateral direction except their independent vertical and rotational displacement. Thus, the 4 wheel parts have 2 DOF each: one allows the wheel to move in vertical direction with regard to the vehicle body, and the other allows the wheel to rotate around the axle.

The generalized coordinates are chosen and denoted with vector $q$:

$$q = \begin{Bmatrix} q_b \\ q_u \end{Bmatrix} = \begin{Bmatrix} [X_{A,b}, Y_{A,b}, Z_{A,b}, \varphi, \phi, \varphi]^T \\ [\theta_{fr}, \theta_{fl}, \theta_{rr}, \theta_{rl}, z_{fr}, z_{fl}, z_{rr}, z_{rl}]^T \end{Bmatrix} \quad (13)$$

The corresponding velocity vector is presented as:

$$\dot{q} = \begin{Bmatrix} \dot{q}_b \\ \dot{q}_u \end{Bmatrix} = \begin{Bmatrix} [\dot{X}_{A,b}, \dot{Y}_{A,b}, \dot{Z}_{A,b}, \dot{\varphi}, \dot{\phi}, \dot{\varphi}]^T \\ [\dot{z}_{fr}, \dot{z}_{fl}, \dot{z}_{rr}, \dot{z}_{rl}, \dot{\theta}_{fr}, \dot{\theta}_{fl}, \dot{\theta}_{rr}, \dot{\theta}_{rl}]^T \end{Bmatrix} \quad (14)$$

The relative position of the wheel center in XY plane of the vehicle reference frame $[\boldsymbol{x}_w, \boldsymbol{y}_w]$ can be denoted as:

$$\begin{Bmatrix} \boldsymbol{x}_w \\ \boldsymbol{y}_w \end{Bmatrix} = \begin{Bmatrix} [l_f, l_f, -l_r, -l_r] \\ \frac{1}{2}[-w_f, w_f, -w_r, w_r] \end{Bmatrix} \quad (15)$$

where $l_i$ and $w_i$ are respectively the $x$ and $y$ positions of each wheel.

The vertical relative position of the unsprung mass $\boldsymbol{z}_w$ in the vehicle reference frame is denoted with its vertical $\boldsymbol{z}_u$, longitudinal $\boldsymbol{x}_w$ and lateral $\boldsymbol{y}_w$ position in the vehicle reference frame:

$$\boldsymbol{z}_w = \boldsymbol{z}_u - (\boldsymbol{y}_w \varphi - \boldsymbol{x}_w \phi + Z_{A,b}) \quad (16)$$

the corresponding velocity of the unsprung mass in the vehicle reference frame can be denoted as:

$$\dot{\boldsymbol{z}}_w = \dot{\boldsymbol{z}}_u - (\boldsymbol{y}_w \dot{\varphi} - \boldsymbol{x}_w \dot{\phi} + \dot{Z}_{A,b}) \quad (17)$$

The motion equations of the 14-DOF vehicle model can be derived based on the Lagrangian dynamics [24]:

$$\begin{cases} \dfrac{d}{dt}\left(\dfrac{\partial T}{\partial \dot{\boldsymbol{q}}_b}\right) - \dfrac{\partial T}{\partial \boldsymbol{q}_b} = \boldsymbol{Q}_b \\ \dfrac{d}{dt}\left(\dfrac{\partial T}{\partial \dot{\boldsymbol{q}}_u}\right) - \dfrac{\partial T}{\partial \boldsymbol{q}_u} = \boldsymbol{Q}_u \end{cases} \quad (18)$$

where $T$ is the kinetic energy of the system, $\boldsymbol{Q}_b$ and $\boldsymbol{Q}_u$ are the generalized forces applied on the sprung mass and unsprung mass respectively.

*1) Kinetic Energy of the Sprung Mass:* In order to fully evaluate the influence of the unsprung mass, the wheels are fixed with the vehicle body in $x_b - y_b$ directions. Based on this consideration, the kinetic energy of the sprung mass can be denoted as:

$$T_b = \frac{1}{2} \boldsymbol{V}_b^T [M_b] \boldsymbol{V}_b + \sum \frac{1}{2} \boldsymbol{V}_{u,i}^T [M_{u,i}] \boldsymbol{V}_{u,i} \quad (19)$$

where the symbols in the above equation will be described in the following paragraphs.

As shown in Fig. 2, there are two reference systems used in this work: the global (inertia) reference system fixed with the ground and the moving reference system fixed on the vehicle body. The origin of the moving frame is located in the mass center of the sprung mass, while the $\boldsymbol{x}_b$, $\boldsymbol{y}_b$ and $\boldsymbol{z}_b$ axles point forward the longitudinal, lateral and vertical direction of motion. The two reference frame are connected with the transformation matrix $[h_{A,b}]$:

$$[h_{A,b}] = \begin{bmatrix} \cos\psi & \sin\psi & 0 & 0 & 0 & 0 \\ -\sin\psi & \cos\psi & 0 & 0 & 0 & 0 \\ 0 & 0 & 1 & 0 & 0 & 0 \\ 0 & 0 & 0 & 1 & 0 & 0 \\ 0 & 0 & 0 & 0 & 1 & 0 \\ 0 & 0 & 0 & 0 & 0 & 1 \end{bmatrix} \quad (20)$$

The velocity of the vehicle body $\boldsymbol{V}_b$ in the moving frame can thus be denoted as:

$$\boldsymbol{V}_b = \begin{bmatrix} V_{xb} \\ V_{yb} \\ V_{zb} \\ \omega_{xb} \\ \omega_{yb} \\ \omega_{zb} \end{bmatrix} = [h_{A,b}] \begin{bmatrix} \dot{X}_{A,b} \\ \dot{Y}_{A,b} \\ \dot{Z}_{A,b} \\ \dot{\varphi} \\ \dot{\phi} \\ \dot{\psi} \end{bmatrix} = [h_{A,b}]\dot{\boldsymbol{q}}_b \quad (21)$$

The velocity $\boldsymbol{V}_{u,i}$ of each unsprung mass is calculated with its relative position in the moving frame and the velocity vector of the vehicle body, which can be denoted as:

$$\boldsymbol{V}_{u,i} = \begin{bmatrix} v_{ux,i} \\ v_{uy,i} \end{bmatrix}$$

$$= \begin{bmatrix} 1 & 0 & 0 & 0 & z_{w,i} & -y_{w,i} \\ 0 & 1 & 0 & -z_{w,i} & 0 & x_{w,i} \end{bmatrix} \boldsymbol{V}_b$$

$$= [h_{b,u,i}][h_{A,b}]\dot{\boldsymbol{q}}_b \quad (22)$$

where $v_{ux,i}$ and $v_{uy,i}$ are, respectively, the longitudinal and lateral velocity of the unsprung mass in the moving frame, while $\boldsymbol{x}_{w,i}$, $\boldsymbol{y}_{w,i}$ and $\boldsymbol{z}_{w,i}$ are the position coordinates of the unsprung mass in the moving frame denoted by Equation (15).

The mass matrix of the sprung mass $[M_b]$ and each unsprung mass $[M_{u,i}]$ are denoted as Equation (23) and Equation (24),

respectively.

$$[M_b] = diag\{m_b, m_b, m_b, J_{xxb}, J_{yyb}, J_{zzb}\} \quad (23)$$

$$[M_{u,i}] = \begin{bmatrix} m_{u,i} & 0 \\ 0 & m_{u,i} \end{bmatrix} \quad (24)$$

where $m_b$ is the mass of the sprung mass, $J_{xxb}$, $J_{yyb}$ and $J_{zzb}$ are, respectively, the inertia of the sprung mass around the $O_b - x_b$, $O_b - y_b$ and $O_b - z_b$ axles.

With the above items, the kinetic energy of the sprung mass can be written in a more compact form:

$$T_b = \frac{1}{2} \boldsymbol{V}_b^T [M_b] \boldsymbol{V}_b + \sum \frac{1}{2} \boldsymbol{V}_{u,i}^T [M_{u,i}] \boldsymbol{V}_{u,i}$$

$$= \frac{1}{2} \dot{\boldsymbol{q}}_b^T [h_{A,b}]^T [M_b][h_{A,b}] \dot{\boldsymbol{q}}_b$$

$$+ \sum \frac{1}{2} \dot{\boldsymbol{q}}_b^T [h_{A,b}]^T [h_{b,u,i}]^T [M_{u,i}][h_{b,u,i}][h_{A,b}] \dot{\boldsymbol{q}}_b$$

$$= \frac{1}{2} \dot{\boldsymbol{q}}_b^T [M_{gb}] \dot{\boldsymbol{q}}_b \quad (25)$$

The generalized mass matrix $[M_{gb}]$ of the sprung mass is a function of the generalized coordinates $Z_{A,b}$, $\psi$ and $z_u$.

*2) Kinetic Energy of the Unsprung Mass:* The kinetic energy of the unsprung mass is composed by the vertical and rotational motion parts, which can be denoted as:

$$T_u = \frac{1}{2} \boldsymbol{\omega}_u^T [J_u] \boldsymbol{\omega}_u + \frac{1}{2} \boldsymbol{V}_{uz}^T [M_u] \boldsymbol{V}_{uz} \quad (26)$$

In this work, the inertia of the propulsion system is taken into account, the angular velocity $\boldsymbol{\omega}_u$ of the propulsion system and the unsprung mass for an independent drive topology can be denoted as:

$$\boldsymbol{\omega}_u = [h_u] \left[ \dot{\theta}_{fr}, \dot{\theta}_{fl}, \dot{\theta}_{rr}, \dot{\theta}_{rl} \right]^T \quad (27)$$

The transformation matrix $[h_u]$ and the inertia matrix $J_u$ are denoted as:

$$\begin{cases} [h_u] = [diag\{1,1,1,1\}; diag\{i_{g,fr}, i_{g,fl}, i_{g,rr}, i_{g,rl}\}] \\ J_u = [diag\{J_{u,fr}, J_{u,fl}, J_{u,rr}, J_{u,rl}, J_{d,fr}, J_{d,fl}, \\ \qquad J_{d,rr}, J_{d,rl}\}] \end{cases} \quad (28)$$

where $J_{u,i}$ and $J_{d,i}$ are respectively the inertia of each wheel and each motor, $i_{g,i}$ is the speed ratio of each transmission.

The vertical velocity matrix $\boldsymbol{V}_{uz}$, and the mass matrix $[M_u]$ of the unsprung mass can be denoted as:

$$\boldsymbol{V}_{uz} = [\dot{z}_{fr}, \dot{z}_{fl}, \dot{z}_{rr}, \dot{z}_{rl}]^T \quad (29)$$

$$[M_u] = diag\{m_{fr}, m_{fl}, m_{rr}, m_{rl}\} \quad (30)$$

With the above terms, the kinetic energy of the unsprung mass is denoted as Equation (31), and it can be derived as a function

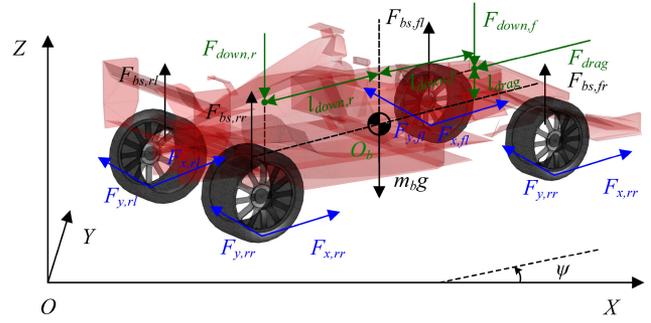

Fig. 3. Vehicle model: forces and torque applied on the sprung mass.

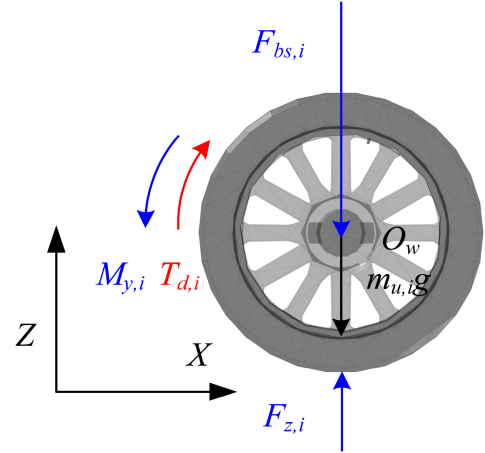

Fig. 4. Vehicle model: forces and torque applied on the unsprung mass.

of the generalized coordinates and generalized mass matrix.

$$T_u = \frac{1}{2} \boldsymbol{\omega}_u^T [J_u] \boldsymbol{\omega}_u + \frac{1}{2} \boldsymbol{V}_{uz}^T [M_u] \boldsymbol{V}_{uz}$$

$$= \frac{1}{2} \begin{bmatrix} \dot{\theta}_{fr} \\ \dot{\theta}_{fl} \\ \dot{\theta}_{rr} \\ \dot{\theta}_{rl} \end{bmatrix}^T [h_u]^T [J_u][h_u] \begin{bmatrix} \dot{\theta}_{fr} \\ \dot{\theta}_{fl} \\ \dot{\theta}_{rr} \\ \dot{\theta}_{rl} \end{bmatrix}$$

$$+ \frac{1}{2} \begin{bmatrix} \dot{z}_{fr} \\ \dot{z}_{fl} \\ \dot{z}_{rr} \\ \dot{z}_{rl} \end{bmatrix}^T [I]^T [M_u][I] \begin{bmatrix} \dot{z}_{fr} \\ \dot{z}_{fl} \\ \dot{z}_{rr} \\ \dot{z}_{rl} \end{bmatrix}$$

$$= \frac{1}{2} \dot{\boldsymbol{q}}_u^T [M_{gu}] \dot{\boldsymbol{q}}_u \quad (31)$$

where $[M_{gu}]$ is the generalized unsprung mass collecting only static values, it is denoted as:

$$[M_{gu}] = \begin{bmatrix} [h_u]^T [J_u][h_u] & \\ & [M_u] \end{bmatrix} \quad (32)$$

*3) Generalized Forces:* The forces and torque acted on the sprung mass and unsprung mass are respectively illustrated as Fig. 3 and Fig. 4. The vertical forces include the gravity of the sprung mass $m_b g$, the four suspension forces $F_{bs,i}$, the front and rear aerodynamic down forces $F_{down,f}$ and $F_{down,r}$, and the anti-roll force $F_{atr,i}$. Forces applied in the longitudinal direction are the four longitudinal tire forces $F_{x,i}$ and the aerodynamic drag force $F_w$. In lateral direction, there are only the four lateral tire forces $F_{y,i}$. The torque applied on the unsprung mass are four driving/braking torque $T_{d,i}$ transmitted by the shafts. The aerodynamic drag and down forces are presented as:

$$\begin{cases} F_{drag} = \frac{1}{2} C_d \rho A V_x^2 \\ F_{down,i} = \frac{1}{2} C_{l,i} \rho A V_x^2 \end{cases} \quad (33)$$

where $C_d$ is the drag coefficient, $\rho$ is the air density, $A$ is the vehicle effective area, $V_x$ is the vehicle longitudinal velocity, $F_{down,i}$ is the aerodynamic downforce, $C_{l,i}$ is the lift coefficient, $i = \{front, rear\}$.

The force matrix including the torque and forces applied on the unsprung mass can be denoted as

$$\boldsymbol{F}_u = \begin{bmatrix} T_{d,fr} + M_{y,fr} - F_{x,fr} z_{fr} \\ T_{d,fl} + M_{y,fl} - F_{x,fl} z_{fl} \\ T_{d,rr} + M_{y,rr} - F_{x,rr} z_{rr} \\ T_{d,rl} + M_{y,rl} - F_{x,fr} z_{rl} \\ F_{atr,fr} - F_{bs,fr} + F_{z,fr} - m_{fr} g \\ -F_{atr,fl} - F_{bs,fl} + F_{z,fl} - m_{fl} g \\ F_{atr,rr} - F_{bs,rr} + F_{z,rr} - m_{rr} g \\ -F_{atr,rl} - F_{bs,rl} + F_{z,rl} - m_{rl} g \end{bmatrix} \quad (34)$$

where $T_{d,i}$ is the driving torque, $M_{y,i}$ is the rolling resistance moment, $F_{atr,i}$ is the anti-roll force, $F_{z,i}$ is the vertical force acted on each tire and $m_i$ is the mass of each wheel.

The force matrix including the forces and torque acted on the sprung mass in each direction is denoted as Equation (35) as shown at the bottom of the page.

The generalized forces can be derived based on the force analysis and virtual work principle, which are presented in the following equations. The similar detail derivation process can be referred to [25].

$$\begin{cases} \boldsymbol{Q}_b = \sum \boldsymbol{F} \frac{\partial \boldsymbol{r}_b}{\partial \boldsymbol{q}_b} = \boldsymbol{F}_b \\ \boldsymbol{Q}_u = \sum \boldsymbol{F} \frac{\partial \boldsymbol{r}_u}{\partial \boldsymbol{q}_u} = \boldsymbol{F}_u \end{cases} \quad (36)$$

*4) Lagrange's Equations:* The generalized motion equations of the rigid sprung and unsprung mass are derived in the form of Equation (37) based on the Lagrangian mechanics and D'Alembert's principle. In this work, we differentiate the generalized mass matrix directly instead of calculating the partial differentials of the system kinetic energy to the time and generalized coordinates separately, $\ddot{\boldsymbol{q}}_b$ and $\ddot{\boldsymbol{q}}_u$ can be denoted directly based on Equation (37) and Equation (38) in this way, which is more efficient for derivation.

$$\frac{d}{dt}\left(\frac{\partial T}{\partial \dot{\boldsymbol{q}}_b}\right) - \frac{\partial T}{\partial \boldsymbol{q}_b} = [M_{gb}]\ddot{\boldsymbol{q}}_b + [\dot{M}_{gb}]\dot{\boldsymbol{q}}_b$$
$$- \frac{1}{2}\dot{\boldsymbol{q}}_b^T\left[\frac{\partial M_{gb}}{\partial \boldsymbol{q}_b}\right]\dot{\boldsymbol{q}}_b = \boldsymbol{Q}_b \quad (37)$$

$$\frac{d}{dt}\left(\frac{\partial T}{\partial \dot{\boldsymbol{q}}_u}\right) - \frac{\partial T}{\partial \boldsymbol{q}_u} = [M_{gu}]\ddot{\boldsymbol{q}}_u + [\dot{M}_{gu}]\dot{\boldsymbol{q}}_u$$
$$- \frac{1}{2}\dot{\boldsymbol{q}}_b^T\left[\frac{\partial M_{gb}}{\partial \boldsymbol{q}_u}\right]\dot{\boldsymbol{q}}_b = \boldsymbol{Q}_u \quad (38)$$

### B. Suspension Model

The suspension model in this section involves the calculation of spring forces, damping forces, anti-roll forces, toe angles, camber angles and the steering angles. The suspension model elaborated below is capable to describe the behavior of both the dependent and independent suspensions.

*1) Spring and Damping Forces:* The suspension force is composed by the spring force and damping force. When the stiffness and damping ratio are constant values, the suspension force can be denoted as Equation (39), the spring force on the spring is denoted as a function of the stiffness $k_{bs,i}$ and the deformation $l_{bs,i}$ of the spring, while the damping force is a function of damping $c_{d,i}$ and velocity $\dot{l}_{d,i}$ of the damper.

$$F_{bss,i} = k_{bs,i}\Delta l_{s,i} + c_{d,i}\dot{l}_{d,i} \quad (39)$$

The deformation of the spring travel $\Delta l_{s,i}$ is a function of wheel jounce $\Delta D_i$, which can be calculated by a transmission

$$\boldsymbol{F}_b = \begin{bmatrix} \sum (F_{x,i}\cos\delta_i - F_{y,i}\sin\delta_i)\cos\psi - \sum (F_{x,i}\sin\delta_i + F_{y,i}\cos\delta_i)\sin\psi - F_w\cos\psi \\ \sum (F_{x,i}\cos\delta_i - F_{y,i}\sin\delta_i)\sin\psi + \sum (F_{x,i}\sin\delta_i + F_{y,i}\cos\delta_i)\cos\psi - F_w\sin\psi \\ \sum F_{bs,i} - mg + F_{down,f} + F_{down,r} \\ -\sum F_{atr,i}y_{u,i} + \sum F_{bs,i}y_{u,i} + \sum (F_{x,i}\sin\delta_i + F_{y,i}\cos\delta_i)Z_{A,b} + \sum T_{d,i}\sin\delta_i \\ \sum F_{down,i}x_{u,i} - \sum F_{bs,i}x_{u,i} - \sum (F_{x,i}\cos\delta_i - F_{y,i}\sin\delta_i)Z_{A,b} - \sum T_{d,i}\cos\delta_i \\ \sum M_{z,i} + \sum (F_{x,i}\sin\delta_i + F_{y,i}\cos\delta_i)x_{w,i} - \sum (F_{x,i}\cos\delta_i - F_{y,i}\sin\delta_i)y_{w,i} \end{bmatrix} \quad (35)$$

ratio $\lambda_{s,i}$,

$$\Delta l_{s,i} = \lambda_{s,i} \Delta D_i \tag{40}$$

where the wheel jounce $\Delta D_i$ is the vertical movement of wheel or axle relative to the vehicle reference frame, which can be defined as:

$$\Delta D_i = z_{w,i} - z_{w0,i} \tag{41}$$

The deformation velocity of the damper can be denoted as:

$$\dot{l}_{s,i} = \lambda_{s,i} \dot{z}_{w,i} \tag{42}$$

where $z_{w,i}$ is the vertical position of the unsprung mass in the moving frame, $z_{w0,i}$ is its initial value.

Finally, the suspension force acted on each wheel can be denoted as:

$$F_{bs,i} = \lambda_{s,i} F_{bss,i} \tag{43}$$

*2) Anti-Roll Forces:* In this work, the anti-roll bars are implemented to reduce the roll displacement of the race car during fast cornering or over road irregularities. The anti-roll bar connects opposite left and right wheels together through a short lever arm linked by a torsion spring. Taking the front axle as an example, the anti-roll force $F_{atr,f}$ is a function of the average jounce $\bar{D}_{l,f}$ and delta jounce $\Delta D_{l,f}$ of the left and right wheels as denoted in Equation (44). The anti-roll forces can also be calculated with the given parameter tables and a 2-D interpretation method.

$$F_{atr,f} = f(\bar{D}_{l,f}, \Delta D_{l,f}) \tag{44}$$

where $\bar{D}_{l,f} = \frac{D_{r,f} + D_{l,f}}{2}$, $\Delta D_{l,f} = D_{r,f} - D_{l,f}$, $D_{r,f}$ and $D_{l,f}$ are respectively the front wheel jounce of the right and left side.

*3) Camber Angles:* The camber angle $\gamma_i$ can be denoted as a function of the wheel jounce and steering wheel angle input by the driver,

$$\gamma_i = f(\Delta D_i, \delta_{driver}) \tag{45}$$

Similarly, the lookup table method based on the 2-D interpolation can be utilized to calculate the camber angle.

*4) Toe Angles:* Toe angle of each wheel $\xi_i$ is also considered in this suspension model, which is denoted as a function of the wheel jounce,

$$\xi_i = f(\Delta D_i) \tag{46}$$

The toe angles can be calculated with the aforementioned 1-D interpolation method.

*5) Steering Angles:* The steering angle of each wheel $\delta_{d,i}$ on the ground is expressed as a function of the wheel jounce and the steering wheel angle input by the driver,

$$\delta_{d,i} = f(\Delta D_i, \delta_{driver}) \tag{47}$$

A 2D interpolation method can be used in this model to obtain the steering angle on the ground $\delta_{d,i}$ with the presented data.

The final steering angle on the ground in the vehicle reference system is:

$$\delta_i = \delta_{d,i} + \xi_i \tag{48}$$

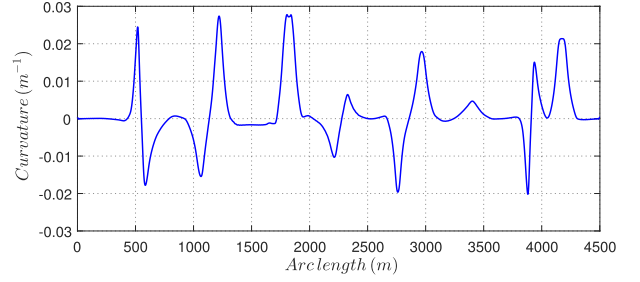

Fig. 5. Curvature of the track.

The tables describing all of the above relationships can be obtained via experiments in bench test for the interpretation method.

### C. Tire Model

In this work, a semi-empirical tire model has been reprogrammed based on the full set of Magic Formula (MF) equations in [26]. To improve the computation efficiency, the MF tire model is programmed in vector format in MATLAB. The longitudinal force $F_{x,i}$, lateral force $F_{y,i}$, overturning moment $M_{x,i}$, rolling resistance moment $M_{y,i}$ and aligning moment $M_{z,i}$ of the tire are calculated with the vertical force $F_{z,i}$, longitudinal slip $\kappa_i$, side slip angle $\alpha_i$ and inclination angle $\gamma_i$, forward velocity $V_x$ as inputs, in this section $i = \{fr, fl, rr, rl\}$ means front right, front left, rear right or rear left.

The $N_p$ pages of forces and torque can be calculated at one time with high computational efficiency by the reprogrammed tire model:

$$[\boldsymbol{F}_x, \boldsymbol{F}_y, \boldsymbol{M}_x, \boldsymbol{M}_y, \boldsymbol{M}_z] = \boldsymbol{f}_{Tire}(\boldsymbol{F}_z, \boldsymbol{\kappa}, \boldsymbol{\alpha}, \boldsymbol{\gamma}, \boldsymbol{V}_x) \tag{49}$$

where the mapped dimensions of the input and output variables are:

$$\boldsymbol{f}_{Tire} : \mathbb{R}^{4 \times N_p} \times \mathbb{R}^{4 \times N_p} \times \mathbb{R}^{4 \times N_p} \times \mathbb{R}^{4 \times N_p}$$
$$\times \mathbb{R}^{1 \times N_p} \to \mathbb{R}^{20 \times N_p}$$

### D. Path Following Model

*1) Curvature of the Track:* The curvature of the Nurburgring circuit can be calculated by Equation (50) with the given X-Y coordinates that can be obtained with GPS, or extracted and converted from the commercial or open source map. The track can be described by its curvature and arc length in a curvilinear coordinate system as it is presented in Fig. 5.

$$C = \frac{dx \cdot ddy - ddx \cdot dy}{\left(\sqrt{dx^2 + dy^2}\right)^3} \tag{50}$$

where $dx$, $ddx$, $dy$, $ddy$ are the first and second order gradients of the X-Y coordinates respectively.

*2) Vehicle Position in Curvilinear Coordinate System:* In curvilinear coordinate system shown in Fig. 6, the vehicular position on the track can be described by its traveled distance $s$ along the reference trajectory, its normal distance to the reference trajectory $n$, and its orientation angle $\theta$ at the current position [27]. The orientation angle can be denoted with the

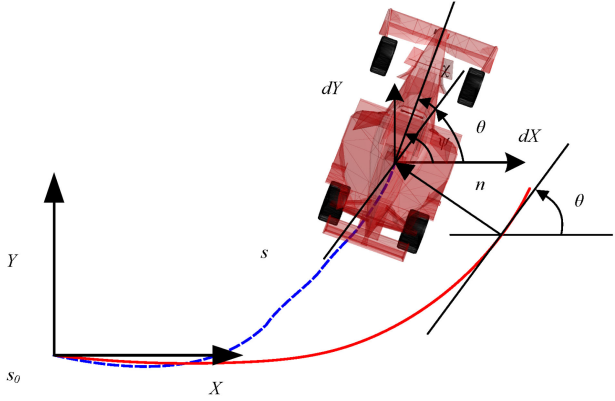

Fig. 6. Vehicle position in curvilinear coordinate system.

yaw angle $\psi$ and the angle between the heading direction and the track $\chi$:

$$\theta = \psi - \chi \tag{51}$$

With Equation (52) and Equation (53), the derivative of the traveled distance $\dot{s}$ can be denoted as Equation (54) with the absolute longitudinal and lateral velocity of the vehicle in global reference frame:

$$\dot{s} - n\dot{\theta} = dX \cos\theta + dY \sin\theta \tag{52}$$

$$\dot{\theta} = \mathcal{C}\dot{s} \tag{53}$$

$$\dot{s} = \frac{dX \cos\theta + dY \sin\theta}{1 - n\mathcal{C}} \tag{54}$$

The derivative of the normal distance $\dot{n}$ can be denoted as:

$$\dot{n} = dY \cos\theta - dX \sin\theta \tag{55}$$

The derivative of $\chi$ can be denoted as:

$$\dot{\chi} = \dot{\psi} - \mathcal{C}\dot{s} \tag{56}$$

### E. Powertrain Mass

In order to evaluate the effect of design parameters of the motor and the transmission on the lap time of the race car, the mass model of the motor and transmission mainly concerning the dependence of the mass and output torque of the powertrain on the design parameters is considered. The total mass of the electric race car is denoted as:

$$m_t = m_b + 4(m_d + m_g) \tag{57}$$

where $m_b$, $m_d$ and $m_g$ are separately the mass of the vehicle body, electric motor and transmission.

The mass of each electric motor and gearbox can be derived as [19]:

$$\begin{cases} m_d = \rho_m \left(\dfrac{P_{\max}}{n_b}\right)^{3/4} \\ m_g = \dfrac{500 T_{in}}{K} \psi \rho_{gc} \rho_g \pi \left(1 + \dfrac{1}{i_g} + i_g + i_g^2\right) \end{cases} \tag{58}$$

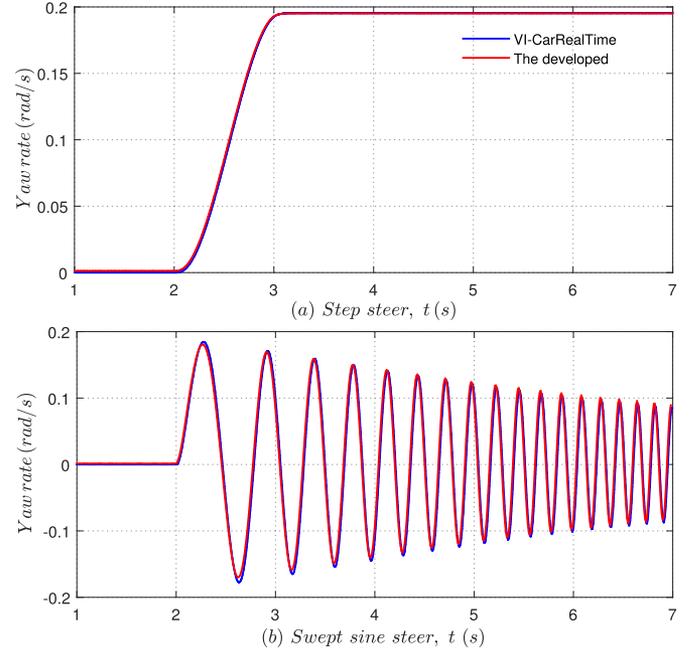

Fig. 7. Comparison of yaw rate in step and swept sine steer manoeuvres.

where $\rho_m$ is the mass factor of the electric motor $P_{\max}$ is the maximum power of the motor (kW), $T_{in}$ is the torque acted on the input shaft of the gearbox (Nm), $K$ is the surface durability factor ($N/mm^2$), $\psi$ is the gear volume fill factor, $\rho_{gc}$ is a correction factor of the gear-box mass model, $\rho_g$ is the mass density of the gear, $i_g$ is the gear ratio.

### F. Model Validation

The 14-DOF model developed in MATLAB is validated in a step and a swept sine steering manoeuvre based on a Formula 3 chassis. The simulation results are validated with a commercial software VI-CarRealTime. The comparisons of the yaw rate are presented respectively in Fig. 7(a) and Fig. 7(b).

## IV. OPTIMAL CONTROL AND OPTIMIZATION SETTINGS

A MATLAB software package for General DYNamic OPTimal control problems abbreviated as *GDYNOPT* is developed to solve the formulated optimal powertrain design and control problems. *GDYNOPT* is developed based on different kinds of transcription methods and differentiation methods, moreover, it is implemented with a novel automatic scaling method and supports parallel computing. The obtained nonlinear optimization problem transcribed by *GDYNOPT* is solved with an open source NLP solver IPOPT [28], the details of *GDYNOPT* which are not the focus of this work will be presented in our future publications. Here we only represent the implemented direct transcription method briefly.

### A. NLP Variables

The continuous state and control variables can be discretized over the whole time interval into $N_n$ nodes with the linear interpretation method. The new independent time samples are obtained by generating a linearly spaced vector on $[t_0, t_f]$ with length $N_n$. The discretized state and control variables at each node together with the static parameters of the continuous optimal design and control problem can be reconstructed as a column vector of NLP variables:

$$\boldsymbol{y} = [\boldsymbol{x}_1, \boldsymbol{u}_1, \boldsymbol{x}_2, \boldsymbol{u}_2, \ldots, \boldsymbol{x}_{N_n}, \boldsymbol{u}_{N_n}, \boldsymbol{t}_f, \boldsymbol{p}]^T \quad (59)$$

where $\boldsymbol{x}_i$ is the row vector of the state variables at node $i$ with dimension $1 \times n_x$, $\boldsymbol{u}_i$ is the row vector of the control variables at node $i$ with dimension $1 \times n_u$, $t_f$ is the final time, $\boldsymbol{p}$ is the row vector of the design variables with dimension $1 \times n_p$. The dimension of $\boldsymbol{y}$ is $N_n(n_u + n_x) + n_{t_f} + n_p$.

### B. Direct Transcription

*1) Defect Constraints Based on Trapezoidal Approach:* The defect constraints $\boldsymbol{\zeta}_k$ of the Trapezoidal method are denoted as:

$$\boldsymbol{\zeta}_k = \boldsymbol{x}_{k+1} - \boldsymbol{x}_k - \frac{1}{2}h(\boldsymbol{f}(\boldsymbol{x}_k, \boldsymbol{u}_k, \boldsymbol{p}) + \boldsymbol{f}(\boldsymbol{x}_{k+1}, \boldsymbol{u}_{k+1}, \boldsymbol{p})) \quad (60)$$

The defect constraints matrix $\boldsymbol{\zeta} \in \mathbb{R}^{(N_n-1) \times n_x}$ can be calculated with a more compact way:

$$\boldsymbol{\zeta} = T_{s2}\boldsymbol{x} - \frac{h}{2}T_{s1}\boldsymbol{f}(\boldsymbol{x}, \boldsymbol{u}, \boldsymbol{p}) \quad (61)$$

where $\boldsymbol{x} \in \mathbb{R}^{N_n \times n_x}$ is the state vector, $\boldsymbol{u} \in \mathbb{R}^{N_n \times n_u}$ is the control variable, $h$ is the time step, $T_{s1}$ and $T_{s2}$ are the transformation matrix, $\boldsymbol{f}(\boldsymbol{x}, \boldsymbol{u}, \boldsymbol{p}) \in \mathbb{R}^{N_n \times n_x}$ is the system dynamics denoted as

$$\boldsymbol{f}(\boldsymbol{x}, \boldsymbol{u}, \boldsymbol{p}) = \begin{bmatrix} \boldsymbol{f}(\boldsymbol{x}_1, \boldsymbol{u}_1, \boldsymbol{p}) \\ \vdots \\ \boldsymbol{f}(\boldsymbol{x}_{N_n}, \boldsymbol{u}_{N_n}, \boldsymbol{p}) \end{bmatrix}$$

$$= \begin{bmatrix} f_1(\boldsymbol{x}_1, \boldsymbol{u}_1, \boldsymbol{p}) & \cdots & f_{n_x}(\boldsymbol{x}_1, \boldsymbol{u}_1, \boldsymbol{p}) \\ \vdots & \ddots & \vdots \\ f_1(\boldsymbol{x}_{N_n}, \boldsymbol{u}_{N_n}, \boldsymbol{p}) & \cdots & f_{n_x}(\boldsymbol{x}_{N_n}, \boldsymbol{u}_{N_n}, \boldsymbol{p}) \end{bmatrix}. \quad (62)$$

*2) Path and Boundary Constraints:* The path constraints $\boldsymbol{g}(\boldsymbol{x}, \boldsymbol{u}, \boldsymbol{p}) \in \mathbb{R}^{N_n \times n_g}$ are the functions of state, control, design and terminal time variables denoted as Equation (63), while the boundary constraints $\boldsymbol{b} \in \mathbb{R}^{1 \times n_b}$ are the functions of initial and final state variables.

$$\boldsymbol{g}(\boldsymbol{x}, \boldsymbol{u}, \boldsymbol{p}) = \begin{bmatrix} \boldsymbol{g}(\boldsymbol{x}_1, \boldsymbol{u}_1, \boldsymbol{p}) \\ \vdots \\ \boldsymbol{g}(\boldsymbol{x}_{N_n}, \boldsymbol{u}_{N_n}, \boldsymbol{p}) \end{bmatrix}$$

$$= \begin{bmatrix} g_1(\boldsymbol{x}_1, \boldsymbol{u}_1, \boldsymbol{p}) & \cdots & g_{n_g}(\boldsymbol{x}_1, \boldsymbol{u}_1, \boldsymbol{p}) \\ \vdots & \ddots & \vdots \\ g_1(\boldsymbol{x}_{N_n}, \boldsymbol{u}_{N_n}, \boldsymbol{p}) & \cdots & g_{n_g}(\boldsymbol{x}_{N_n}, \boldsymbol{u}_{N_n}, \boldsymbol{p}) \end{bmatrix} \quad (63)$$

TABLE I
SIMULATION PARAMETERS OF THE POWERTRAIN OPTIMIZATION

| Prameters | Values |
|---|---|
| Mass factor of the motor $\rho_m$ | 0.2845 |
| Gear volume fill factor $\psi$ | 0.7 |
| Mass density of the gear $\rho_g$ | 78550 $kg/m^3$ |
| Sprung mass with powertrain $m_v$ | 490 $kg$ |
| Distance of mass center to front axle $l_f$ | 1.595 $m$ |
| Distance of mass center to rear axle $l_r$ | 1.135 $m$ |
| Front wheel track $w_f$ | 1.585 $m$ |
| Rear wheel track $w_r$ | 1.535 $m$ |
| Maximum power $P_{max}$ | 165 $kW$ |
| Aerodynamic drag force location $h_g$ | 0.18 $m$ |
| Surface durability factor $K$ | 8920000 $N/mm^2$ |
| Mass factor of the gearbox $\rho_{gc}$ | 3.1136 |

After the presented constraints are obtained, the NLP constraints are ready to be constructed as:

$$\boldsymbol{c}(\boldsymbol{y}) = [\boldsymbol{\zeta}_{1,1}, \ldots, \boldsymbol{\zeta}_{1,n_x}, \ldots, \boldsymbol{\zeta}_{(N_n-1),1}, \ldots, \boldsymbol{\zeta}_{(N_n-1),n_x},$$
$$\boldsymbol{g}_{1,1}, \ldots, \boldsymbol{g}_{1,n_g}, \ldots, \boldsymbol{g}_{N_n,1}, \ldots, \boldsymbol{g}_{N_n,n_g},$$
$$\boldsymbol{b}_{1,1}, \ldots, \boldsymbol{b}_{1,n_b}]^T, \ \boldsymbol{c}(\boldsymbol{y}) \in \mathbb{R}^{((N_n-1)n_x + N_n n_g + n_b) \times 1} \quad (64)$$

Both the lower and upper bounds of the defect constraints are zeros, while the bounds of path and boundary constraints are set by the user. All the constraints are reconstructed consistently into the required format of the NLP solver. The obtained NLP aims to minimize the objective function by finding the vector $\boldsymbol{y}$:

$$\min F(\boldsymbol{y}) \quad (65)$$

subject to the bounds and constraints:

$$\begin{cases} \boldsymbol{y}_{\min} \leq \boldsymbol{y} \leq \boldsymbol{y}_{\max} \\ \boldsymbol{c}_{\min} \leq \boldsymbol{c}(\boldsymbol{y}) \leq \boldsymbol{c}_{\max} \end{cases} \quad (66)$$

### C. Simulation Parameters and Settings

The researched electric race car is based on the chassis of a Formula 3 race car, the simulation parameters of which are presented in Table I.

For the NLP solver, the constraints violation is set as $1e-7$ which is acceptable for all the constraints, the desired converging tolerance of the formulated problem is set as $1e-3$. The Jacobians of the defect and path constraints are evaluated with the central differential method, while the Hessian matrices are approximated with the inbuilt limited-memory BFGS approach of the NLP solver [28]. The computation of the optimal design and control problems were performed with MATLAB 2015a on a Linux OS based cluster with Intel Xeon CPU X5355 @2.66 GHz. The computation efficiency is approximately improved by 3 times with parallel computation using 8 CPUs on each node.

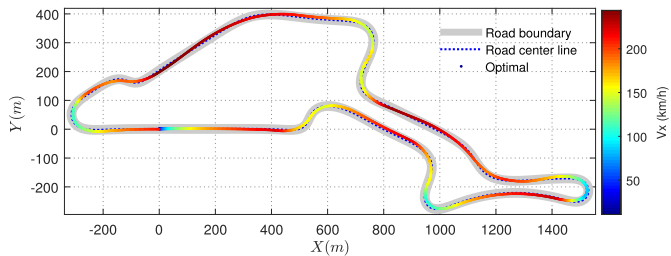

Fig. 8. The optimal racing line of 4-IWD electric race car.

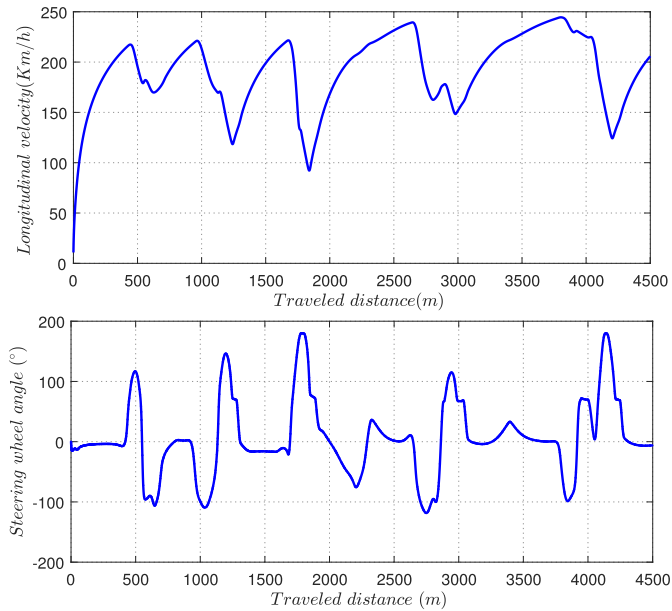

Fig. 9. The optimal longitudinal velocity profile and steering wheel angle.

## V. OPTIMAL POWERTRAIN DESIGN

### A. Design Variables

The electric powertrain design in this work is implemented with four uniform motors and single speed transmissions. The base speed $n_b$ and the constant power speed ratio (CPSR) $\beta$ are selected as the design parameters of the motor. The envelope curve of the torque-speed characteristics of the motor can be obtained with the two design parameters and the given maximum power of the motor $P_{\max}$. The gear speed ratio $i_g$ is selected as the design variable of each transmission. The powertrain design vector of the 4-IWD electric race car is:

$$\boldsymbol{p} = [n_b, \beta, i_g], \ n_p = 3 \tag{67}$$

### B. Design Results

The optimal racing line of the case with on-board motor is demonstrated in Fig. 8, while Fig. 9 presents the longitudinal velocity profile and optimal steering wheel angle respectively.

The optimal powertrain design results with the on-board and in-wheel motors are demonstrated in Table II. The interesting thing is that, the electric race car with 4 in-wheel motors can achieve a better lap time performance. The improvement is 0.281 s on the test track compared with the on-board motor configuration. The influence of the unsprung mass on the lap time performance is a little different from the fact that addressed by most engineers and researchers. More research work should be conducted to have a deep insight.

TABLE II
OPTIMAL POWERTRAIN DESIGN RESULTS

| Prameters | $n_b$ | $\beta$ | $i_g$ | $m_p$ (kg) | $t_f$ |
|---|---|---|---|---|---|
| On-board motor | 6195 | 3.92 | 7.86 | 46.55 | 88.418 |
| In-wheel motor | 6208 | 4.02 | 8.60 | 50.25 | 88.227 |

TABLE III
SENSITIVITY OF THE LAP TIME TO THE UNSPRUNG MASS

| Unsprung mass (kg) | Sprung mass (kg) | Inertia ($kgm^2$) | Lap time (s) |
|---|---|---|---|
| [11, 11, 13, 13] | 490 | 0.5 | 88.238 |
| [15, 15, 17, 17] | 474 | 0.5 | 88.178 |
| [19, 19, 21, 21] | 458 | 0.5 | 88.112 |
| [23, 23, 25, 25] | 442 | 0.5 | 88.054 |

TABLE IV
SENSITIVITY OF THE LAP TIME TO THE UNSPRUNG ROTATIONAL INERTIA

| Unsprung mass (kg) | Sprung mass (kg) | Inertia ($kgm^2$) | Lap time (s) |
|---|---|---|---|
| [11, 11, 13, 13] | 490 | 0.5 | 88.223 |
| [11, 11, 13, 13] | 490 | 0.75 | 88.651 |
| [11, 11, 13, 13] | 490 | 1 | 88.803 |
| [11, 11, 13, 13] | 490 | 1.25 | 88.925 |

*1) Sensitivity to the Unsprung Mass:* In order to analyze the influence of increasing the unsprung mass of the 4-IWD electric race car on the lap time performance, the powertrain parameters are fixed with the obtained ones in last section. Optimal control will be applied to find only the control parameters that minimizing the lap time of the 4-IWD electric race car with different unsprung mass. The obtained results are illustrated in Table III. In these three cases, different values of sprung mass are moved to the unsprung mass, while the rotational inertia is a uniform constant value which can be realized by reasonably changing the shape of the unsprung masses. It is illustrated in Table III that when more unsprung mass is moved to the sprung mass, the lap time is decreased, which is different from the everybody knows conclusions. The underlying reason will be analyzed and presented in next section.

*2) Sensitivity to the Unsprung Inertia:* In this subsection, the unsprung mass is fixed and the lap time of race cars with different unsprung rotational inertia are analyzed with optimal control. As we can see from Table IV, the lap time performance of the race car is very sensitive to the unsprung rotational inertia, which will get worse with the increasing of the rotational inertia.

## VI. OPTIMAL CHASSIS DESIGN

The section will explore the influences of mass center, anti-roll bar and suspension on the lap time of the on-board 4-IWD electric race car with the optimal control theory.

### A. Optimization of the Mass Center

The mass center of the chassis is known to have significant influences on the lap time. Considering this, the distance from the mass center to front axle of the chassis $l_f$ is added into the design parameters in this work, the design parameters of the 4-IWD electric race car become into $[\beta, n_b, i_g, l_f]$. For the on-board electric race car, the newly obtained optimal design parameters are $[3.16, 7604, 8.35, 1.344]$, compared with the data in Table I, the mass center of the race car is moved forward by 0.251 m. Fig. 10 demonstrates that the normal loads acted on the front tires increase while the ones acted on the rear wheels decreases and the normal loads of the race car with the optimized mass center distribute more equally. The obtained lap time is $t_f = 86.880$, which is improved by 1.538 s compared with the obtained results in Table II. For the in-wheel electric race car, the newly obtained optimal design parameters are $[3.17, 7226, 8.28, 1.313]$ and the achieved lap time is $t_f = 86.850$. We can see that after optimization of the mass center of the race car the obtained lap time of the on-board and in-wheel race car cases are very similar now. According to the results in III, the mass distribution condition may get worse and result in poorer lap time performance when only the unsprung mass is moved to sprung.

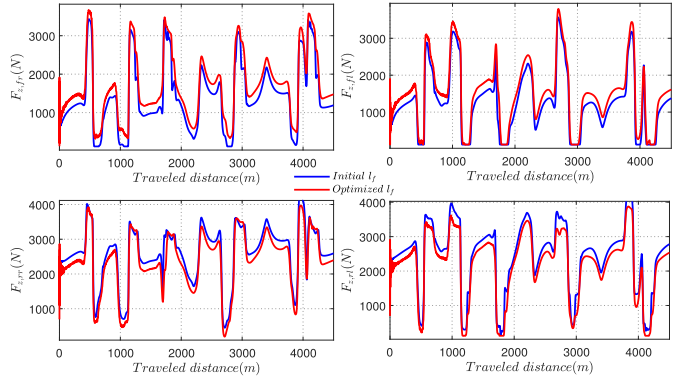

Fig. 10. The normal loads acted on the four tires.

### B. Using Two Different Pairs of Motors

The average normal loads and longitudinal propulsion forces distributed on the front and rear tires are not equal in the test maneuver, which is the fact in almost all the maneouvers due to different amount of accelerating and braking operations. However, currently, the four motors and transmissions are implemented with the same physical parameters in almost all the developed four wheel driving electric vehicles, which is also the assumption of the above optimization in this work. In this case, the power of the motors mounted on the front axle is a kind of surplus while the rear axle suffers power deficit in the tested maneouver. The achieved solution above still has not taken full advantage of the power of each motor which limits further improvement of the lap time, especially considering the fact that the propulsion power of the race car is limited by the race event.

Based on the above analysis, we propose to utilize two different pairs of motors and transmissions mounted on the front and rear axles respectively. Each pair has two same motors and transmissions but different with the other pair. The updated design parameters of the on-board 4-IWD electric race car are given as:

$$\boldsymbol{p} = [\beta_f, n_{bf}, i_{gf}, \beta_r, n_{br}, i_{gr}, P_r, l_f], n_p = 8 \quad (68)$$

where $\beta_f, n_{bf}, i_{gf}$ and $\beta_r, n_{br}, i_{gr}$ are respectively the CPSR ratio, base speed of the motors and the speed ratio of the transmissions mounted on the front and rear axles, $P_r$ (kW) is the maximum power of the rear motors, the power of the front motors is $P_f = P_{\max} - P_r$.

The achieved optimal design parameters of the 4-IWD electric race car with two different pairs of motors are: $\boldsymbol{p} = [3.19, 7186, 8.91, 3.63, 6330, 7.72, 120, 1.343]$. The achieved solution utilizes two motors with small power (22.5 kW) driving the front wheels and two bigger motors (60 kW) driving the rear wheels. The distance from the mass center to the front axle is 1.343 m. The obtained lap time is 86.462 s which is further improved by 0.418 s compared with the previous design.

### C. Optimization of the Anti-Roll Bar

The anti-roll bars implemented to reduce the roll displacement of the race car when performing a fast cornering or running on a bumpy road are also known to have influence on the lap time. Two more proportional coefficients of the original front $k_{f,atr}$ and rear $k_{r,atr}$ anti-roll bars' parameters are put into the design vector:

$$\boldsymbol{p} = [\beta_f, n_{bf}, i_{gf}, \beta_r, n_{br}, i_{gr}, P_r, l_f, k_{f,atr}, k_{r,atr}], n_p = 10 \quad (69)$$

The corresponding achieved optimal design result is $\boldsymbol{p} = [2.62, 8916, 9.30, 4.33, 5763, 7.14, 125.3, 1.508, 0, 10]$, which tends to use a very stiff anti-roll bar for the rear axle but do not use anti-roll bar for the front axle. The obtained lap time is reduced to 85.710 s.

### D. Optimization of the Suspension

The two main functions of a suspension are to maintain grip by keeping the tires in contact with the road and to provide comfort to the passengers, especially when a dip or a bump appears. The tuning of the suspension of a race car is also of great importance in improving the lap time. The design variables are finally augmented with four more parameters which are the stiffness value of the front ($k_{bs,f}, N/m$) and rear ($k_{bs,r}, N/m$) springs, proportional coefficients of the original front $c_f$ and rear $c_r$ damping ratios.

$$\boldsymbol{p} = [\beta_f, n_{bf}, i_{gf}, \beta_r, n_{br}, i_{gr}, P_r, l_f,$$
$$k_{bs,f}, k_{bs,r}, c_f, c_r, k_{f,atr}, k_{r,atr}], n_p = 14 \quad (70)$$

The obtained optimal design parameters are: $\boldsymbol{p} = [2.52, 9409, 9.56, 4.30, 5803, 7.31, 132, 1.690, 2.54 \times 10^5, 1 \times 10^6, 0.66, 9.09, 0.08, 10]$, and the achieved lap time is $t_f = 85.225$, which is further improved by 0.485 s. Many people think that relatively soft springs are important for grip performance by intuitive reasoning. The optimization result is a little 'shock', since the optimal solution tends to utilize very stiff spring for front suspension and 'super stiff' spring for the rear suspension. The underlying reasons may be: 1) Stiffer rear springs will make the car more likely to oversteer; 2) It can help to overcome some additional weight loading due to weight transfer, stiffer springs will resist harder against compression and help to prevent the traction loss; 3) Soft rear spring may result in better handling performances, however, it may be not worth to make this compromise since aerodynamics does not prefer and there is rare advantage in improving ride comfort for the race car drivers. Besides, the damping ratios of the front dampers are increased while the ones of the rear are opposite in the optimized suspension.

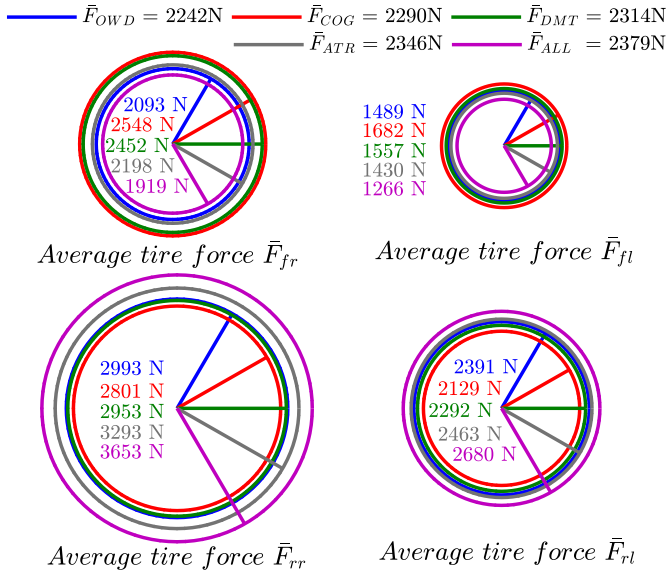

Fig. 11. Average tire forces of different optimization cases.

### E. Comparisons of the Five Optimization Cases

In order to compare the improvements of the successive optimizations, the average tire forces of the four tires and the total average grip forces in different optimization cases are respectively calculated by Equation (71) and Equation (72).

$$\bar{F}_i = \frac{1}{N_n} \sum_{k=1}^{N_n} \sqrt{F_{tx,i}^2(k) + F_{ty,i}^2(k)} \quad (71)$$

$$\bar{F}_j = \frac{1}{N_n} \sum_{k=1}^{N_n} \sqrt{\sum_{i=1}^{4}(F_{tx,i}^2(k) + F_{ty,i}^2(k))/4} \quad (72)$$

where $j = \{OWD, COG, DMT, ATR, ALL\}$, $OWD$ means the case of optimizing only the powertrain parameters, the subsequent symbols mean the optimizations with augmented parameters of mass center, different pairs of motors, anti-roll bar and suspensions in order.

The comparisons of the five optimization cases are demonstrated in Fig. 11.

Several findings can be presented: 1) The average grip forces of the right side are more than the left side, which is caused by the more left steering operations; 2) The average grip forces of the front tires are less than the rear, which is caused by the accelerating operations and the mass distribution; 3) After optimizing the mass center, the forces acted on the front axle increases while the ones on the rear are decreased. The grip force of each tire tends to distribute more equally since the four implemented motors are uniform and the optimization tries to take full advantages of all the motors. The total average grip force is increased by 48 N, the lap time is 86.850 s which is improved by 0.432 s; 4) After using different pairs of motors, the total average grip force and the lap time are respectively further increased by 24 N and 86.462 s, the lap time improvement is 0.388 s; 5) After introducing the optimization of the anti-roll bar parameters, there are more grip forces from the rear tires but less from the front. The lap time becomes 85.710 s with a significant improvement of 0.752 s; 6) When the optimization parameters are augmented with the ones of the suspension, the tire forces on the rear increases and the ones on the front decreases further. The total average grip force is increased by 6.1% in comparison with the case of optimizing only the powertrain parameters, the achieved lap time is 85.225 s which is improved again by 0.485 s. The step by step optimizations can be considered as the process of increasing the grip forces of the tires to the limit.

## VII. CONCLUSIONS

In this work, the optimal design and control problems of a 4-IWD electric race car are investigated based on a developed 14-DOF vehicle model and an optimal control software package *GDYNOPT*. The in-wheel 4-IWD electric race car with increased unsprung mass is found to obtain very similar or even better lap time performance than the on-board case with uniform rotational inertia. The lap time performance is very sensitive to the rotational inertia, which can be reduced by adjusting the shape of the wheels and propulsion system. The optimization of the mass center can improve the lap time performance by distributing the normal load more appropriately in order to take the utmost use of all tires. The parameters of anti-roll bars also have significant influences on the lap time, using very stiff rear anti-roll bar has improved the lap time a lot in this work. The optimized suspensions using very stiff front springs and super stiff rear springs can help to reduce the lap time further. The obtained optimal design parameters in this work can be utilized to design or tune the components of the race car, while the obtained optimal trajectory and longitudinal velocity profile can serve as the references for training the race car drivers. The design of a race car in practice is far more complicated which is generally full of compromise, the future work may take into account more details and experimental work will be done when a accessible platform is built. In addition, the developed vehicle model and software package can also be utilized in other kinds of design and control problems with different optimization purposes.

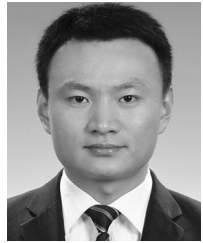

**Huilong Yu** (M'17) received the M.Sc. degree in mechanical engineering from the Beijing Institute of Technology, Beijing, China, and the Ph.D. degree in mechanical engineering from the Politecnico di Milano, Milano, Italy, in 2013 and 2017, respectively. He is currently a Research Fellow of advanced vehicle engineering with the University of Waterloo, Waterloo, ON, Canada. His research interests include vehicle dynamics, optimal control, closed-loop control, and energy management problems of conventional, electric, hybrid electric, and autonomous vehicles.

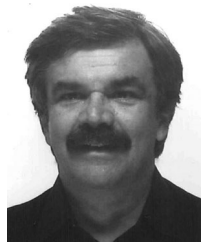

**Federico Cheli** received the Graduate degree in mechanical engineering from the Politecnico di Milano, Milano, Italy, in 1981. He is currently a Full Professor with the Department of Mechanics, Politecnico di Milano. From 1992 to 2000, he was an Associate Professor with the Faculty of Industrial Engineering, Politecnico di Milano. He is the Co-founder of the E_CO spin-off and an author of more than 380 publications in international journals or presented at international conferences. His scientific activity concerns research on vehicle performance, handling and comfort problems, active control, ADAS, and electric and autonomous vehicles. He is a member of the editorial board of the *International Journal of Vehicle Performance* and *International Journal of Vehicle Systems Modeling and Testing*.

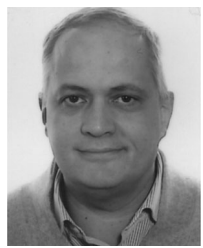

**Francesco Castelli-Dezza** (M'94) received the M.Sc. and Ph.D. degrees in electrical engineering from the Politecnico di Milano, Milano, Italy, in 1986 and 1990, respectively. He is currently a Full Professor with the Department of Mechanical Engineering, Politecnico di Milano. His research interests include studies on dynamic behavior of electrical machines, electrical drives control and design, and power electronics for energy flow management. He is a member of the IEEE Power Electronics Society.